\documentclass{aoscm}
\input{epsf}
\newcommand{\IR}{\mathbb R}

\newcommand{\Il}{{\bf 1}}

\newtheorem{theorem}{Theorem} 

\newtheorem{remark}[theorem]{Remark}

\let\epsilon=\varepsilon

\let\hat=\widehat

\let\tilde=\widetilde 
\renewcommand{\theequation}{\arabic{equation}}

\begin{document}

\SPECFNSYMBOL{1}{2}{}{}{}{}{}{}{} 
\AOSMAKETITLE

\AOSyr{2000}
\AOSvol{00}
\AOSno{00}
\AOSpp{000--000}
\AOSReceived{Received }
\AOSAMS{Primary 62G07, 62C20; Secondary 60G70, 41A25}
\AOSKeywords{Wavelets, Thresholding, Minimax}
\AOStitle{Stein Estimation for Infinitely Divisible Laws}
\AOSauthor{R. Averkamp\thanks{Research supported in part
by the NSF Grant DMS-9632032, while this 
author was visiting the Southeast Applied Analysis Center, 
School of Mathematics, Georgia Institute of Technology, 
Atlanta, GA 30332, USA.}
and C. Houdr\'e\thanks{Research supported in part
by NSF and NSA Grants.}}
\AOSaffil{Freiburg University
and Universit\'e Paris XII and Georgia Institute of Technology}
\AOSlrh{R. AVERKAMP and C. HOUDR\'E}
\AOSrrh{STEIN ESTIMATION}
\AOSAbstract{Unbiased risk estimation, \`a la Stein, 
is studied for infinitely
divisible laws with finite second moment}

\maketitle
\BACKTONORMALFOOTNOTE{5}

Let us start by briefly  
recalling the framework and results of Stein (\cite{S}): Let $X_i$, 
$i=1,\ldots,n$, be iid $N(0,\sigma^2)$ random variables and let 
$g=(g_1,\ldots,g_n):\IR^n \longrightarrow \IR^n$, be ``weakly 
differentiable."  Then for all 
$\theta \in \IR^n$,  
\begin{eqnarray} \label{stein}
E\|X+\theta +g(X+\theta)-\theta\|^2_2= n \sigma^2 &+&E\|g(X+\theta)\|^2_2\\
&+& 2 \sigma^2 E\sum_{i=1}^n 
\frac{\partial }{\partial x_i} g_i(X+\theta),\nonumber
\end{eqnarray}
where $\|\cdot\|_2$ is the Euclidean norm.  
Thus the risk of the estimator $x+g(x)$ can be estimated unbiasedly 
by $\displaystyle n \sigma^2 +g(x)^2+2 \sigma^2 \sum_{i=1}^n 
 \frac{\partial g_i}{\partial x_i}(x)$.
This estimate is useful only if the variance of 
the risk estimate is small compared to the actual risk.
This is especially the case if $g_i$ only depends on $X_i$,
since then the strong law of large numbers kicks in.  
For normal random variables the 
existence of the above estimates is based on the 
identity
\[
\int_\IR g^\prime (x) e^{-x^2/2} dx =\int_\IR x g(x) e^{-x^2/2} dx.
\]  
We obtain below a corresponding identity for
infinitely divisible random variables with finite variance, replacing
$g^\prime$ by $K(g)$, where $K$ is an operator commuting
with translations.

Let $f$ be a density on $\IR$, with mean $0$ and variance $\sigma^2$  
(for simplicity of notation 
we concentrate on the univariate case, but see Remark~\ref{multi}).  
Let $d(x)=x+g(x)$ be an estimator in the location model induced by $f$.
Let $F=\{f*\delta_\theta: \theta\in \IR\}$, 
while $L^2(F)$ and $L^1(F)$ have their canonical
meaning.
We want to estimate the risk of $d$ unbiasedly:
\begin{eqnarray*}
\lefteqn{\int_\IR (d(x+\theta)-\theta)^2f(x)dx}\\
&&=\int_\IR g(x+\theta)^2 f(x)dx +\int_\IR x^2f(x)dx + 
2\int_\IR x g(x+\theta) f(x)dx.
\end{eqnarray*}
In the above right hand side, the first summand can be 
estimated unbiasedly, 
the second is a constant,
so we just need to find a function $h\in L^1(F)$ such that
\begin{equation}
\label{h-equation}
\int_\IR h(x+\theta)f(x)dx=\int_\IR x g(x+\theta) f(x)dx.
\end{equation}
If $g$ is a polynomial the right-hand side of 
(\ref{h-equation}) is itself a 
polynomial in $\theta$.  It is then well-known that 
there exists an $h$ 
satisfying (\ref{h-equation}).  But if $g$ is the soft--thresholding 
operator, i.e.,  
$g(x) =T^S_\lambda(x)=(|x|-\lambda)^+{\rm sgn}(x)$, then $g$ does 
not even have a power series expansion.  
Moreover, $h$ does not have to be unique.  Indeed, 
$h+q$ is also a solution for 
any function $q$ such that $q*f=0$ 
(if $\hat f$, the Fourier transform of $f$, 
has zeros such $q$ might exists).

Hence, let us assume that $\hat f$ does not have zeros.
By computing the generalized Fourier transform of both sides of
\begin{equation}
\int_\IR g(-x+\theta) (-x) f(-x)dx=\int_\IR h(-x+\theta)f(-x)dx,
\end{equation}
we get:
\begin{equation}
\label{K-hat-equation}
\hat{g}(w) \hat{f}^\prime (-w) =i\hat{h}(w) \hat{f}(-w).
\end{equation}
This identity shows that if 
$\hat g$ converges to $0$ fast enough, e.g., if 
$\hat{g}$ has  compact support, then there exists an $h$ such that
(\ref{h-equation}) holds.  
Since $\hat f$ does not vanish, $h$ is uniquely 
determined.
Hence the set 
\begin{eqnarray*}
\lefteqn{U_f:=}\\
&&\hspace{-0.5cm}
\left\{\!g\in L^2(F): \exists h \in L^1(F) , 
\int_\IR\!h(x+\theta)\!f(x)dx=\int_\IR\!\!x 
g(x+\theta)\!f(x)dx,\!\,\forall \theta \in \IR \!\right\}, 
\end{eqnarray*}
is a vector space and clearly there is a unique linear map 
$K_f:U_f \longrightarrow L^1(f)$ with 
\[ \int_\IR K_f(g)(x+\theta)f(x)dx=\int_\IR g(x+\theta)xf(x)dx.\]

Let us note some properties of $K_f$:
\begin{theorem}
\label{K-properties}
Let $f, f_1, f_2$ be densities with finite second moment, and let 
$K_f$, $K_{f_1}$ and $K_{f_2}$ be well defined.  
Then for all $b \in \IR$ and $g\in U_f$, respectively 
$g \in U_{f_1*f_2}$:
\begin{enumerate}
\item $K_f(g(\cdot+b))=K_f(g)(\cdot+b)$
\item $K_{f*\delta_b}(g) =K_f(g)+b g(\cdot)$
\item $K_{f_1*f_2} (g)=K_{f_1}(g)+K_{f_2}(g)$
\item $K_{b f(\cdot b)}(g) =K_f(g(\cdot/b)) (\cdot b)/b$, for $b>0$.
\end{enumerate} 
\end{theorem}
\proof{Proof}
1. $\int_\IR g(x+\theta +b)x f(x)dx = 
\int_\IR K_f(g)(x+\theta +b)f(x)dx$ and hence 
$K_f(g(\cdot+b))=K_f(g)(\cdot+b)$.\\
2.
\begin{eqnarray*}
\int_\IR x g(x+ \theta) f(x-b)dx&=& 
\int_\IR (K_f(g)(x+\theta+b)+b g(x+\theta+b))f(x)dx\\
&=& \int_\IR (K_f(g)(x+\theta) +b g(x+\theta)) f(x-b)dx.
\end{eqnarray*}
3. let $h_1$, $h_2$ be such that 
$\int_\IR g(x+\theta)x f_i(x)dx = \int_\IR h_i(x+\theta)f(x)dx$, $i=1,2$,
then 
\begin{eqnarray*}
\lefteqn{\int_\IR (h_1+h_2)(z+\theta) (f_1*f_2)(z)dz}\\
 &=& \int_\IR \int_\IR (h_1(x+y+\theta) +h_2(x+y+\theta)) f_1(x)f_2(y)dxdy\\
&=&\!\int_\IR\!\int_\IR\!h_1(x+y+\theta) f_1(x) dx f_2(y) dy 
+\int_\IR\!\int_\IR\!h_2(y+x+\theta) f_2(y) dy f_1(x) dx\\
&=& \int_\IR g(z+\theta) z (f_1*f_2)(z)dz.
\end{eqnarray*}
4. 
\begin{eqnarray*}
\int_\IR  g(x+\theta) x (b f(b x))dx &=& \int_\IR g(x/b+\theta) x f(x)/b dx\\
&=& \int_\IR g((x+b \theta)/b) x/b f(x)dx\\
&=& \int_\IR K_f(g(\cdot/b)) (x+b\theta)/b f(x) dx\\
&=& \int_\IR K_f(g(\cdot/b)) ((x+\theta)b)/b (bf(x b))dx,
\end{eqnarray*}
and thus $K_{b f(\cdot b)}(g) =K_f(g(\cdot/b)) (\cdot b)/b$.
\endproof  

Note that the third property presented above is  
very useful for wavelet analysis, since the law of the 
noise in a wavelet coefficient is a weighted convolution of 
the noise in the original data.

For the normal distribution with unit variance the operator $K$ is
defined by $K(g)=g^\prime$, i.e. $K$ is the differentiation operator.  
In general $K$ is quite complicated to compute, however from
(\ref{K-hat-equation}) we see that formally 
\[\widehat{K_f(g)}(w)=\hat{g}(w)\hat{f}^\prime (-w)/(\hat{f}(-w)i).\]
This suggest that $h$ can be computed by a convolution of the estimator
and of a function or measure, which is the inverse Fourier 
transform of $\hat{f}^\prime(-w)/(\hat{f}(-w)i)$.  Let us try to further 
formalize this claim.  Assume $K_f(g):=K_f*g$, where $K_f\in L^1(\IR)$
and $\hat K_f= \hat f^\prime(-\cdot)/(\hat f(-\cdot)i)$. 
If $g \in L^\infty(\IR)$, then $K_f*g$ does what it is supposed to do:
\begin{eqnarray*}
\int_\IR (K_f*g) (x+\theta) f(x) dx &=& 
\int_\IR \int_\IR K_f(x-t) g(t+\theta)dt f(x) dx\\
&=& \int_\IR \int_\IR K_f(x-t) f(x)dx g(t+\theta)dt \\
&=& \int_\IR \int_\IR K_f(-(t-x)) f(x)dx g(t+\theta)dt\\
&=& \int_\IR  (K_f(-\cdot)*f(\cdot))(t) g(t+\theta)dt\\
&=& \int_\IR t f(t) g(t+\theta)dt,
\end{eqnarray*}
where the last equality follows from the construction of $K$:
\[\hat{K}(-\cdot) \hat{f}=
{\hat{f}^\prime\over i}=\widehat{(f(\cdot)\mbox{id})}.\]

If $K_f$ is known, then it can still be a problem to compute $K_f(g)$, since 
$K_f(g)$ is not necessarily as simple as $g^\prime$.  
If $g=\sum_i g_i$ and the $K_f(g_i)$ are easy to compute, then 
we can compute $K_f(g)$ since $K_f$ is linear.  
For example we can take 
$g^+_\lambda(x)=(x-\lambda)^+$ and $g_\lambda^-(x):=(x-\lambda)^-$ 
as simple building blocks for functions.  
Note that $K_f(g_\lambda^+)(x)=K_f(g_0^+)(x-\lambda)$ and 
$K_f(g_0^-(x))=\sigma^2-K_f(g_0^+)$ if $\int_\IR x f(x)dx=0$ and 
$\int_\IR x^2 f(x)dx=\sigma^2$.
For example the soft thresholding estimator $T^S_\lambda$ 
as well as $T^M_\lambda$ given by 
$T^M_\lambda(x):= x \Il_{\{|x|\ge\lambda\}} + 2 (|x|-\lambda/2)_+ 
{\rm sgn(x)} \Il_{\{|x|<\lambda\}}$, have the 
following decompositions:
\begin{equation}
\label{soft-comp}
T^S_\lambda(x) = x -g^+_0(x) + g^+_\lambda(x)-g_0^-(x)+g^-_{-\lambda}(x),
\end{equation}
\[
 T^M_\lambda(x) =x - g_0^+(x) +2 g^+_{\lambda/2}(x) -g_\lambda^+(x)
                   - g^-_0(x)+2 g^-_{-\lambda/2}(x) -g^-_{-\lambda}(x) .\]
For a further example, assume that 
$g :\IR^+\longrightarrow \IR$, is twice continuously 
differentiable with $g(0)=0$, 
then $g(x)=g^\prime(0^+)x^+ +\int_0^\infty (x-y)^+ g^{\prime\prime}(y)dy$.

Another simple example is provided by compound Poisson distributions.  Indeed, 
let $F$ be a compound Poisson distribution with Fourier transform 
$\exp(\lambda(\Psi(w)-1))$,
where $\Psi$ is the characteristic function of the density $f$.
Then 
$$\hat K_F= \frac{\lambda\Psi^\prime(-w)}{i},$$ and thus 
$K_F(x)= -\lambda f(-x) x $,
i.e. $K_F(g)=K_F*g$.  Since compound distributions are building blocks for 
infinitely divisible ones, we have:  

\begin{theorem}
Let $f$ be an infinitely divisible density with finite second moment, i.e., let 
\[ \hat f(t)= \exp\left(ibt + 
\int_\IR \left(\frac{\exp(ixt)-1-ixt}{x^2}\right) M(dx) 
\right),\]
where $M$ is a finite positive measure. Let $M(\{0\})=0$, let $b=0$ and let 
$g$ be Lipschitz.  Then 
\[K(g)(t):=\int_\IR \frac{g(t+x)-g(t)}{x} M(dx),\]
is a well defined real valued function, which is moreover 
bounded and continuous.  Furthermore, 
\[  \int_\IR K(g) (x+\theta) f(x)dx =\int_\IR xg(x+\theta) f(x) dx.\]
\end{theorem}

\proof{Proof}
It is clear that $K(g)$ is well defined, bounded and continuous.
If $g$ has compact support then 
$\int_\IR |g(x+y)-g(x)|/|y| M(dy)$ and $K(g)$ are in $L^1(\IR)$  and 
\begin{eqnarray*}
\widehat{K(g)}(t)&=& 
\int_\IR \int_\IR \frac{g(y+x)-g(y)}{x}M(dx) \exp(i t y) dy\\
&=& \int_\IR\int_\IR \frac{g(y+x)-g(y)}{x} \exp( i t y) dy M(dx)\\
&=& \hat g(t)  \int_\IR \frac{\exp(-i x t) -1}{x} M(dx).
\end{eqnarray*}
Since $\int_\IR (\exp(-i x t) -1)/x M(dx) =\hat f^\prime(-t)/(\hat f(-t)i)$,
the Fourier transforms of $\int_\IR K(g)(x+\theta) f(x) dx$ and 
$ \int_\IR g(x+\theta) x f(x) dx$ are equal 
and thus these two terms are themselves equal 
for all $\theta \in \IR$.

If $g$ is Lipschitz but does not have compact support, 
then let $g_n(x):=(1-|x|/n)_+ g(x)$.  Then the Lipschitz 
constants of the $g_n$ form a bounded set.  
Clearly, $g_n$ and $K(g_n)$ respectively 
converge pointwise, respectively to $g$ and $K(g)$, 
and moreover $\|K(g_n)\|_\infty$ is bounded. Hence
$\lim_{n\longrightarrow \infty} \int_\IR g_n(x+\theta) x f(x) dx=
\int_\IR g(x+\theta) x f(x) dx$
and $\lim_{n \longrightarrow \infty} \int_\IR K(g_n)(x+\theta) f(x) dx = 
\int_\IR K(g)(x+\theta) f(x) dx$, for all $\theta$.
Thus 
\[  \int_\IR K(g) (x+\theta) f(x)dx =\int_\IR g(x+\theta)x f(x) dx.\]
\endproof 

\begin{remark}\rm
The assumption $b = 0$, is not serious, $b$ is a location parameter
of the density and thus we can use Theorem \ref{K-properties}.
The condition $M(\{0\})=0$ is not restrictive either.
If $M(\{0\})=\sigma^2$, then the distribution is the convolution of 
a centered normal distribution with variance $\sigma^2$ and 
of an infinitely divisible 
distribution with L\'{e}vy measure without atom at the origin.
Again we can use Theorem \ref{K-properties} for this situation.
Hence in the general case we obtain:
\[ K(g)(t):= b g(t)+ M(\{0\})  g^\prime(t) +\int_{\IR\backslash \{0\}} 
\frac{g(t+x)-g(t)}{x} M(dx).\]
We also note here that although of little interest to us since 
we are dealing with mean square errors, the operator $K$ could as well be defined just under 
a finite first moment assumption on $X$ (in which case, the requirement on $M$ will change too).   
\end{remark}

\begin{remark}\rm
Let $f$ be a density with mean zero and variance $\sigma^2$.
Without loss of generality let also $K_f(1)=0$.
Let $f_n= *_{i=1}^n \sqrt{n} f(\cdot \sqrt{n})$.
By the central limit theorem $f_n$ converges in distribution to a
normal density. So one would expect that $K_{f_n}$ converges in some
sense to $\sigma^2 d/dx$. 
Assume that $K_f(g)(x)= \int_\IR (g(x+y)-g(x))/y M(dy)$.  Note that if 
$K_f(g)=Q*g$, where $Q$ is a measure, then with the notation 
$Q^-(A):=Q(-A)$,
\begin{eqnarray*} 
(Q*g)(x)&=&\int_\IR g(x-y)-g(x)Q(dy)\\
&=& \int_\IR \frac{g(x+y)-g(x)}{y} y Q^-(dy),
\end{eqnarray*}
where the first equality holds since $K_f(1)=0$, i.e. $\int_\IR 1 Q(dx)=0$.
Since $\int_\IR x (x+\theta) f(x)dx = \int_\IR x^2 f(x)dx$,
taking $g(x)=x$ gives $M(\IR)=K(x)$ and 
$\int_{\IR} K(x)f(x)dx=\int_\IR x^2 f(x)dx=\sigma^2$.
As we already know
$ K_{f_n} (g)(x)$ \linebreak $ = K_f(g(\cdot/\sqrt{n})) 
(x \sqrt{n}) /\sqrt{n}$.
Using the form of $K_f$, we now have 
\[K_{f_n}(g)(x) =\int_\IR \frac{g(x+y/\sqrt{n})-g(x)}{y/\sqrt{n}} M (dy).\]
Thus, if $g$ is Lipschitz and differentiable, 
$\lim_{n \longrightarrow \infty} K_{f_n}(g)(x)=\sigma^2 g^\prime(x)$.
\end{remark}

\proof{Examples}
1. Let $f(x)=\exp(-\sqrt{2}|x|)/ \sqrt{2}$ be 
the variance normalized Laplace 
density.
It is easy to see that,
$ \hat{f}(w)=2/(2+w^2)$.
Thus
\[{\hat{f}^\prime(w)\over i\hat{f}(w)}={2 i w \over 2+w^2}\,,\]
and
\[ \hat{K}(w)={-2 i w \over 2+w^2}=- i w \hat{f}(w)= 
   \widehat{f^\prime}. \]
Thus $K(x)=-\exp(-\sqrt{2} |x|) \mbox{sgn}(x)\in L^1(\IR)$.
Tedious but simple computations yield 
\[K*x_+= \left\{ 
\begin{array}{ccl}
{\displaystyle \frac{\exp(\sqrt{2}x)}{2}}& :& x\le 0\\
{\displaystyle 1- \frac{\exp(-\sqrt{2}x)}{2}} &:& x>0
\end{array}
\right\}=:h(x).\]
Using (\ref{soft-comp}) we obtain
\begin{eqnarray*}
K(T^S_\lambda(x)-x)&=&-h(x)-(1-h(x))+h(x-\lambda)+(1-h(x+\lambda))\\
&=& h(x-\lambda)-h(x+\lambda).
\end{eqnarray*}
Combining these results, we see that for $X$ with a Laplace
distribution,
\[E(T^S_\lambda(X+\theta)-\theta)^2=1+E \min((X+\theta)^2,\lambda^2) + 
         2(h(X+\theta-\lambda)-h(X+\theta+\lambda)).\]

\medskip\noindent
2. Let $f_t(x)= \exp(-x)x^{t-1} /\Gamma(t)\Il_{\IR^+}(x)$ 
be the density of the Gamma
distribution. Since the mean of this distribution is $t$ we want to compute
$K_{f_t*\delta_{-t}}$. 
Then by Feller \cite[p.~567]{Fe}, $\log(\hat f_t(x))=t \int_0^\infty 
(\exp(iyx) -1)/y \exp(-y)dy$ and thus 
$(\log \hat f_t)^\prime(x)= it \int_0^\infty \exp(iyx) \exp(-y)dy$.
Thus $K_{f_t} (g)=Q*g$ where $Q \in L^1(\IR)$ and 
\[\hat Q(x)=t\int_0^\infty \exp(-iyx) \exp(-y)dy = t\int_{-\infty}^0 \exp(ixy)
\exp(y)dy.\]
Hence $K_{f_t*\delta_{-t}}(g)(x)=t\int_{-\infty}^0 \exp(y) g(x-y)dy-tg(x)$.   
Now, symmetrizing $f_t$ (to have a zero mean density) gives 
$\tilde f_t(x)= \exp(-|x|)|x|^{t-1} /2\Gamma(t)\Il_{\IR}(x)$ from which the corresponding 
$K_{\tilde f_t}$ follows.  
\\[0.3cm]
3. Another example is the cosine hyperbolic 
density, $f(x) =1/\cosh(\pi x/2)$,
again \cite[p.~567]{Fe} 
\[\log(\hat f(x))= \int_{\IR} \frac{\exp(ixy)-1-iyx}{y^2}\,\,
\frac{y}{\exp(y)-\exp(-y)} dy\]
and thus 
\[K_f(g)(x)=
\int_\IR \frac{g(x+y)-g(x)}{y}\,\, \frac{y}{\exp(y)-\exp(-y)}dy. \]
\endproof

\begin{center}
\begin{figure}[ht]
\epsfbox{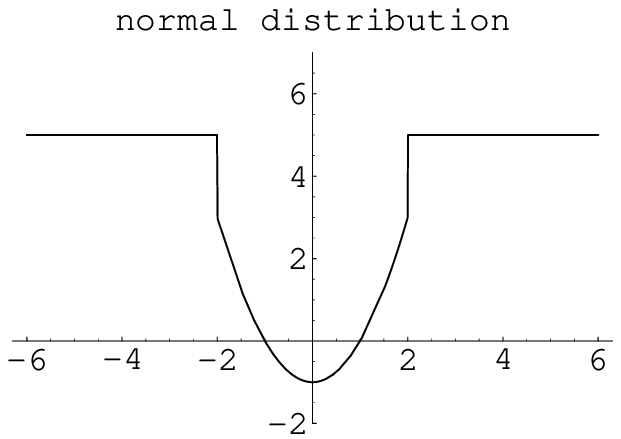}
\epsfbox{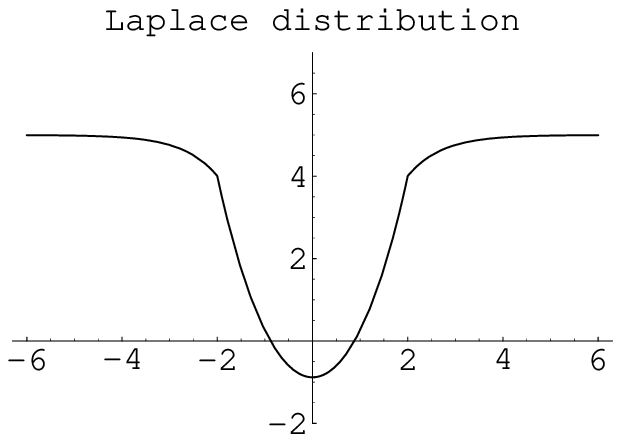}
\epsfbox{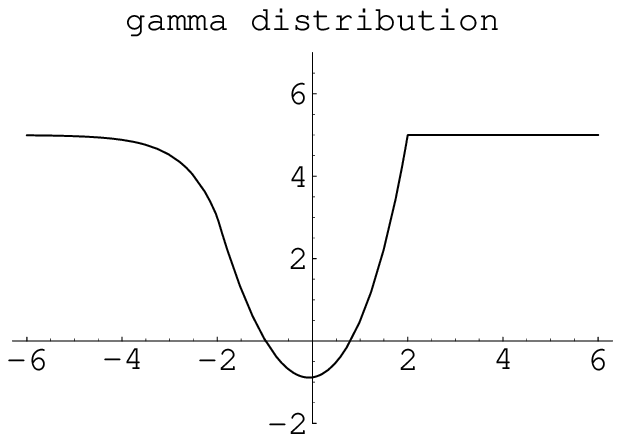}
\epsfbox{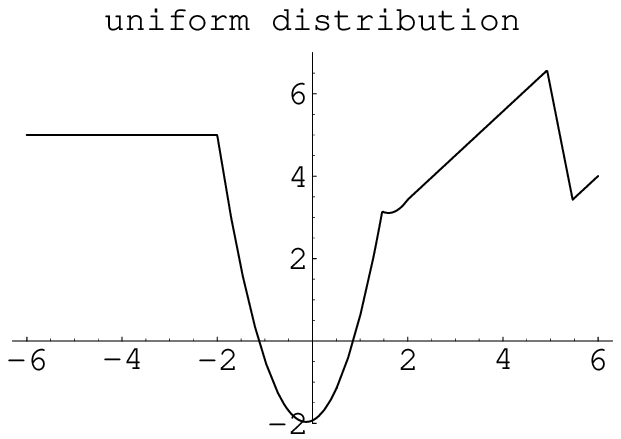}
\caption{The unbiased risk estimators for soft thresholding}
\label{unbiased-risk-estimators}
\end{figure}
\end{center}

The examples presented above are 
infinitely divisible distributions and so $K$ has a nice form.   
Let us consider a case which is not:  The uniform distribution 
with density $2^{-1}\Il_{(-1,1)}$.  
Assume $g:[-1,1]\longrightarrow \IR$ 
and $2^{-1}\int_{-1}^1 g(x)dx=0$. If $\bar g$ is the
2--periodic extension of $g$ on $\IR$ then $\bar 2^{-1}g * \Il_{(-1,1)}=0$.
Thus unbiased risk estimators are not uniquely determined. 
Let 
\[r(\theta)=  2^{-1}\int_{-1}^1 {x(x+\theta)_+} dx = 
\cases{0 & : $\theta \le -1$\cr
\frac{1}{6} + \frac{\theta}{4} - 
\frac{\theta^3}{12} & : $\theta \in (-1,1)$\cr
\frac{1}{3} & : $\theta\ge 1$}.\]
After some tries one finds that with
\[ h(x)= \cases{ 0 & : $x \le 0$\cr
             -\frac{(x-[x/2]2)(x-[x/2]2 -2)}{2} &: $x \ge 0$ }.\]
($h$ is the 2--periodic extension of $-x(x-2)/2$ defined on $[0,2]$ to 
$\IR^+$.)
$2^{-1}\int_{-1}^1 h(x+\theta)dx =r(\theta)$.
So with the help of (\ref{soft-comp}) we can now compute 
an unbiased risk estimator for soft thresholding.
Figure~\ref{unbiased-risk-estimators} shows the unbiased risk estimators 
for soft thresholding
with threshold 2 for the normal distribution, 
the Laplace distribution, the gamma 
distribution with $t=2$ and the uniform distribution. The
distributions were transformed to have unit variance and mean zero.

\newpage

\begin{remark}\label{multi}\rm
As we have seen with (\ref{stein}), for normal random variables,
unbiased risk estimation is possible for multivariate means, even if
the estimators for the coordinates are not independent. 
This is also possible
for other types of distributions, one has to apply the operator $K$ 
coordinatewise.
Let $X_i$, $i=1,\ldots,n$  be random variables, $X_i$ has distribution
$F_i$ and $EX_1=0$, $EX_1^2=\sigma_1^2$. 
Assume that an operator $K_1$ exists such that 
$EX_1g(X_1+\theta_1)=EK_1(g)(X_1+\theta_1)$ for some $g$.
If $g:\IR^n \longrightarrow \IR$ and $\theta \in \IR^n$ then
$E (X_1+g(X+\theta)-\theta_1)^2=\sigma_1^2
+E g(X+\theta)^2+2 E X_1g(X+\theta)$.
Then, under the proper conditions on $g$,
\begin{eqnarray*}
\lefteqn{EX_1 g(X+\theta)}\\
  &=& \int_{\IR^{n-1}}\int_\IR x_1g(x_1+\theta_1,\ldots,x_n+\theta_n) F_1(dx) 
    \otimes_{i=2}^n F_i (d (x_2,\ldots,x_n))\\
&=& \int_{\IR^{n-1}} 
\int_\IR K_1(g(\cdot,x_2+\theta_2,\ldots,x_n+\theta_n))(x_1+\theta_1) 
F_1(dx)\\
   && \hspace{5cm} \hfill \otimes_{i=2}^n F_i(d (x_2,\ldots,x_n)).
\end{eqnarray*}
Thus $E(X_1+g(X+\theta)-\theta_1)^2$ $=\sigma_1^2+E g(X+\theta)^2+$ 
$2 EK_1(g(\cdot,X_2+\theta_2,\ldots,X_n+\theta_n)(X_1+\theta_1)$.
\end{remark}

\begin{remark}\rm
As the reader might have guessed by now, the motivation for the present 
paper comes from thresholding methods in wavelet denoising 
(see \cite{DJKP}).  
In a function space approach to denoising, the thresholds depend on 
the sample size $n$, on the Besov space to which the target 
functions belong to and also on the Besov norm of these targets.
In practice it is often not known which threshold is appropriate since 
the function space to which the signal belongs as well as 
the value of its norm are unknown.  To bypass this problem, 
Donoho and Johnstone developed a procedure called SureShrink where 
thresholds are chosen 
automatically (see \cite{DJ95}). Their method, based on Stein's 
unbiased risk estimate is as follows:  
for each level (except the highest levels) 
in the noisy wavelet transform, the largest threshold 
(smaller than $\sqrt{2 \log n}$) which
minimizes the unbiased risk estimate is chosen.  
For soft thresholding finding this
minimum is simple and takes $O(n \log n)$ time. 

As noticed in (\cite{AH1}), the central 
limit theorem works fast for 
wavelet coefficients, so it is reasonable to
apply the normal adaptive results to the general non Gaussian framework. 
However, it is also of interest to understand the scope of SureShrink 
beyond the normal framework.  To do so, we needed to find 
unbiased risk estimates for other types of distributions.  
This is what we did here for infinitely divisible and related noise. 
We could then potentially use the corresponding 
thresholds found in \cite{AH1} and \cite{AH2}.   
 
\end{remark}

\Line{\AOSaddress{
R.\ Averkamp\\
Institut f\"ur mathematische Stochastik\\
Freiburg University\\
Eckerstra{\ss}e 1\\
79104 Freiburg, Germany\\
}\hfill
\AOSaddress{
C.\ Houdr\'{e}\\
Laboratoire d' Analyse et de\\ 
Math\'ematiques Appliqu\'ees\\
CNRS UMR 8050\\
Universit\'e Paris XII\\
94010 Cr\'eteil Cedex, France\\
and\\
School of Mathematics\\
Georgia Institute of Technology\\
Atlanta, GA 30332, USA}}


\begin{thebibliography}{99}

\bibitem{AH1}  Averkamp, R. and Houdr\'e, C.  
Wavelet Thresholding for non necessarily 
Gaussian Noise: Idealism.  
{\it Annals of Statistics \bf 31}, 110-151 (2003).  

\bibitem{AH2}  Averkamp, R. and Houdr\'e, C.  
Wavelet Thresholding for non necessarily 
Gaussian Noise: Functionality. 


\bibitem{DJ95}  Donoho, D.L. and Johnstone, I.M., 
  Adapting to Unknown Smoothness via Wavelet Shrinkage.
{\it J. Amer. Statist. Assoc. \bf 90}, 1200-1224 (1995).

\bibitem{DJKP}
Donoho, D.L., Johnstone, I.M., Kerkyacharian, G., and Picard, D.,
Wavelet Shrinkage: Asymptotia?.
{\it J. Roy. Statist. Soc. Ser. B, \bf 57}, 301--369 (1995).

\bibitem{Fe} Feller, W., 
{\it An Introduction to Probability  Theory and its Applications},
Vol. II, John Wiley \& Sons (1966).


\bibitem{S} Stein, C., Estimation of the mean of a 
multivariate normal distribution. 
{\it The Annals of Statistics \bf 9}, 
1135-1151 (1981).  


\end{thebibliography}
\end{document}